\documentstyle[12pt]{ article}
\begin{document}

\newtheorem{theorem}{Theorem}[section]
\newtheorem{prop}[theorem]{Theorem}
\newenvironment{mtheorem}{\begin{prop}\rm}{\end{prop}}
\newtheorem{coroll}[theorem]{Corollary}
\newenvironment{mcorollary}{\begin{coroll}\rm}{\end{coroll}}
\newtheorem{lemm}[theorem]{Lemma}
\newenvironment{mlemma}{\begin{lemm}\rm}{\end{lemm}}
\newtheorem{anoprop}[theorem]{Proposition}
\newenvironment{mproposition}{\begin{anoprop}\rm}{\end{anoprop}}
\newenvironment{mproof}{\begin{trivlist}\item[]{\em
Proof: }}{\hfill$\Box$\end{trivlist}}
\newenvironment{mdefinition}{\begin{trivlist}\item[]
{\em Definition.}}{\hfill$\Box$\end{trivlist}}
\newtheorem{eg}{\rm\sl \uppercase{Example}}[section]
\newenvironment{example}{\begin{eg}\rm}{\hfill$\Box$\end{eg}}

\newcommand{\ca}{{\cal A}}
\newcommand{\cb}{{\cal B}}
\newcommand{\cc}{{\cal C}}
\newcommand{\cd}{{\cal D}}
\newcommand{\cf}{{\cal F}}
\newcommand{\cg}{{\cal G}}
\newcommand{\ch}{{\cal H}}
\newcommand{\ck}{{\cal K}}
\newcommand{\cl}{{\cal L}}
\newcommand{\cn}{{\cal N}}
\newcommand{\cm}{{\cal M}}
\newcommand{\co}{{\cal O}}
\newcommand{\cs}{{\cal S}}
\newcommand{\ct}{{\cal T}}

\newcommand{\csa}{{$C^*$-algebra}}
\newcommand{\nsa}{{non-self-adjoint}}
\newcommand{\cssa}{{$C^*$-subalgebra}}

\def \IR{\hbox{{\rm I}\kern-.2em\hbox{{\rm R}}}}
\def \iR{\hbox{{\sevenrm I\kern-.2em\hbox{\sevenrm R}}}}
\def \IN{\hbox{{\rm I}\kern-.2em\hbox{\rm N}}}
\def \IC{\hbox{{\rm I}\kern-.6em\hbox{\bf C}}}
\def \IQ{\hbox{{\rm I}\kern-.6em\hbox{\bf Q}}}
\def \ZZ{\hbox{{\rm Z}\kern-.4em\hbox{\rm Z}}}

\newcommand{\ga}{{\alpha}}
\newcommand{\gb}{{\beta}}
\newcommand{\gd}{{\delta}}
\newcommand{\gA}{{\Alpha}}
\newcommand{\gB}{{\Beta}}
\newcommand{\gG}{{\Gamma}}
\newcommand{\gD}{{\Delta}}
\newcommand{\gs}{{\sigma}}
\newcommand{\gS}{{\Sigma}}

\newcommand{\ot}{\otimes}
\newcommand{\op}{\oplus}
\newcommand{\pr}{\prime}

\newcommand{\bmi}{\bf $i$}
\newcommand{\bmj}{\bf $j$}

\begin{center}
\rm

{\Large \bf On the Banach space isomorphism type of \\ AF C*-algebras
and their triangular \\ subalgebras}
\vspace{.3in}

{\Large S. C. Power}

\vspace{.15in}
\rm

\it Department of Mathematics\\
Lancaster University\\
England LA1 4YF
\rm

\end{center}

\vspace{.3in}
\begin{center}
\bf ABSTRACT
\end{center}

\rm It is shown that all the approximately finite dimensional
C*-algebras which are not of Type I are isomorphic as Banach spaces.
This generalises the matroid case given previously by Arazy.
Analogous results are obtained for various families of triangular
subalgebras of AF C*-algebras.
In addition the classification of various continua of
Type I AF C*-algebras is discussed.
\vspace{.3in}
\rm
\newpage

A C*-algebra is approximately finite (or AF) if it is
a closed union of
finite-dimensional C*-subalgebras. Those AF C*-algebras for which
the  finite-dimensional subalgebras can be taken to be full matrix
algebras are known as matroid C*-algebras. J. Arazy \cite{ara} has shown that
with the exception of the algebra $\cal{K}$ of compact operators, which
is the unique matroid algebra with separable dual, all
(infinite-dimensional) matroid
C*-algebras are isomorphic as Banach spaces. We generalise this
by showing that all the AF C*-algebras which are not of Type I are
isomorphic.
In particular we see that a simple AF
C*-algebra is either isomorphic to $\cal{K}$ or to the Fermion algebra
$F$, the unital matroid C*-algebra $\displaystyle{\lim_\to}\
M_{2^k}$.
We also comment on the classification of AF C*-algebras of Type I.

A similar analysis is given  for various distinguished triangular subalgebras
of AF C*-algebras and it is this non-self-adjoint context
that motivated the present study.
In particular the refinement limit algebras
$\rho - \displaystyle{\lim_\to}\  T_{n_k}$ are all isomorphic to the
model algebra $\ct_{2^\infty} = \rho -\displaystyle{\lim_\to}\  T_{2^k}$, and
the standard limit algebras are all isomorphic to the model algebra
$\cs_{2^\infty} = \sigma -\displaystyle{\lim_\to}\ T_{2^k}$. From this
we are able to deduce that all the (proper) alternation algebras are
isomorphic. It
seems probable that $\ct_{2^\infty}$ and $\cs_{2^\infty}$ are not isomorphic.
Some evidence in support of this is that there are
no ``natural'' complemented contractive injections
$F \rightarrow
\cs_{2^\infty}$ and for this reason the methods of this paper cannot be
applied.
In fact it seems likely that
non-self-adjoint subalgebras of AF
C*-algebras provide an interesting diversity of Banach spaces.
This is to be expected in view of the analogies that exist
between triangular algebras and function spaces. Also, in an associated
nonselfadjoint context, Arias \cite{ari} has recently shown that
there are nonisomorphic nest subalgebras (perhaps uncountably many)
in the trace class.

The  proofs below are straightforward and self-contained.
In particular we
do not require the C*-algebraic classification of the matroid algebras or
the AF C*-algebras (given in Dixmier \cite{dix} and Elliott \cite{ell}
respectively).  The
key step in the self-adjoint case is to show that if $B$ is an AF
C*-algebra whose Bratteli diagram has a particular property, which we
call the Fermion property, then for any other AF C*-algebra $A$
there exists a complemented linear contractive injection $\gamma : A
\rightarrow B$.  This injection is obtained by constructing an infinite
commuting diagram whilst
at the same time constructing a commuting system of left inverses for
the partial injections. This ensures that in the limit the closed
space $\gamma (A)$ is complemented in terms of a completely
contractive linear map.
This form of explicit construction
lends itself to natural modifications to deal with triangular
subalgebras of AF C*-algebras.

As this paper was being prepared Simon Wasserman pointed out to the author
that the isomorphism of non Type I AF algebras has also been obtained
in a different way by Kirchberg as a consequence of his work on
nuclearly embeddable C*-algebras and exact C*-algebras.  In fact,
Kirchberg deduces that all separable nuclear non Type I C*-algebras
are isomorphic. The rather deep result that leads to this isomorphism
is that separable unital nuclear C*-algebras are unitally completely
isometrically embeddable in the Fermion algebra. The proof of this, as
well as the simplified proof given by Wasserman \cite{was}, is quite involved
with C*-algebra theory and contrasts markedly with our approach
for the simpler case of AF C*-algebras.

I would like to thank Jonathan Arazy for telling me about his
classification paper during the Durham Symposium on The
Geometry of Banach Spaces
and Operator Algebras in 1992. I would also like to thank Gordon Blower,
Graham Jameson and Przemek Wojtaszczyk for some helpful conversations.
\newpage

\section{AF algebras with The Fermion property}

The Fermion algebra (or CAR algebra) is the unital matroid C*-algebra
$F = \displaystyle{\lim_\to}\  M_{2^k}$ with Bratteli diagram

\begin{center}
\setlength{\unitlength}{0.0125in}%
\begin{picture}(5,154)(80,660)
\thicklines
\put( 85,800){\line( 0,-1){ 40}}
\put( 85,740){\line( 0,-1){ 40}}
\put( 80,740){\line( 0,-1){ 40}}
\put( 80,800){\line( 0,-1){ 40}}
\multiput( 80,675)(0.00000,-7.50000){3}{\makebox(0.4444,0.6667){\tenrm .}}
\put( 80,685){\makebox(0,0)[lb]{\raisebox{0pt}[0pt][0pt]{\twlrm 8}}}
\put( 80,805){\makebox(0,0)[lb]{\raisebox{0pt}[0pt][0pt]{\twlrm 2}}}
\put( 80,745){\makebox(0,0)[lb]{\raisebox{0pt}[0pt][0pt]{\twlrm 4}}}
\end{picture}
\end{center}

\noindent A general AF C*-algebra, with presentation $A =
\displaystyle{\lim_\to}\  A_k$, has an associated Bratteli diagram in
which, similarly, multiple edges between two vertices indicate, through
their multiplicity, the multiplicity of the partial embedding between
the summands associated with the vertices. We shall say that
the Bratteli diagram has the {\em Fermion property} if there is
a sequence of vertices $v_1, v_2, \dots$ associated with summands of
$A_{n_1}, A_{n_2}, \dots$, respectively,
with $n_k$ an increasing sequence, such that
there is more than one vertical path between each consecutive pair $v_i,
v_{i+1}$. That is, there is a sequence of {\em partial} embeddings
$A_{n_k} \rightarrow A_{n_{k+1}}$, with nonzero compositions, each of
multiplicity at least two. Thus, the Pascal triangle Bratteli diagram
has the Fermion property, whereas the following diagram does
not:

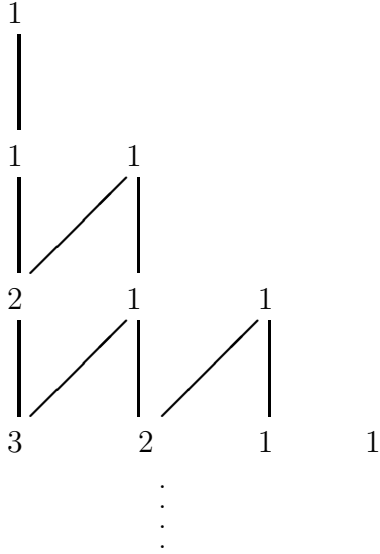
\begin{center}
\setlength{\unitlength}{0.0125in}%
\begin{picture}(150,229)(75,585)
\thicklines
\put(130,740){\line( 0,-1){ 40}}
\put( 85,640){\line( 1, 1){ 40}}
\put(130,680){\line( 0,-1){ 40}}
\put(140,640){\line( 1, 1){ 40}}
\put(185,680){\line( 0,-1){ 40}}
\put( 80,800){\line( 0,-1){ 40}}
\put( 80,740){\line( 0,-1){ 40}}
\put( 80,680){\line( 0,-1){ 40}}
\multiput(140,610)(0.00000,-8.33333){4}{\makebox(0.4444,0.6667){\tenrm .}}
\put( 85,700){\line( 1, 1){ 40}}
\put( 75,745){\makebox(0,0)[lb]{\raisebox{0pt}[0pt][0pt]{\twlrm 1}}}
\put(225,625){\makebox(0,0)[lb]{\raisebox{0pt}[0pt][0pt]{\twlrm 1}}}
\put( 75,805){\makebox(0,0)[lb]{\raisebox{0pt}[0pt][0pt]{\twlrm 1}}}
\put( 75,685){\makebox(0,0)[lb]{\raisebox{0pt}[0pt][0pt]{\twlrm 2}}}
\put(125,745){\makebox(0,0)[lb]{\raisebox{0pt}[0pt][0pt]{\twlrm 1}}}
\put(125,685){\makebox(0,0)[lb]{\raisebox{0pt}[0pt][0pt]{\twlrm 1}}}
\put(180,685){\makebox(0,0)[lb]{\raisebox{0pt}[0pt][0pt]{\twlrm 1}}}
\put(180,625){\makebox(0,0)[lb]{\raisebox{0pt}[0pt][0pt]{\twlrm 1}}}
\put( 75,625){\makebox(0,0)[lb]{\raisebox{0pt}[0pt][0pt]{\twlrm 3}}}
\put(130,625){\makebox(0,0)[lb]{\raisebox{0pt}[0pt][0pt]{\twlrm 2}}}
\end{picture}
\end{center}

We now show that two AF C*-algebras whose Bratteli diagrams have the
Fermion property are linearly homeomorphic. The next lemma, which is
analogous to Lemma 2.11 of \cite{ara}, is the key result required in the proof.

Let $A$ and $B$ be finite-dimensional C*-algebras with chosen matrix
unit systems $\{e_{ij} : (i,j) \in I\}$ and $\{f_{ij}: (i,j) \in J\}$
respectively. Assume that $I, J$ are block diagonal subsets of $\{1,
\dots, m\}^2$ for some $m$. Let $a = (a_{ij})$ belong to $A$. A linear
map $\gamma: A \rightarrow B$ is said to be of {\it compression type}
with respect to these systems if $\gamma$ is a (block diagonal) direct
sum of maps of the form $\alpha \circ \beta$ where $\beta : A
\rightarrow M_n$ is given by

\[
\beta ((a_{ij})) = (a_{k_s,k_t})_{s,t=1}^n,
\]

\noindent where $\{k_1, \dots, k_n\}^2 \subseteq I$, and where $\alpha : M_n
\rightarrow B$ is a multiplicity one algebra injection of the form

\[
\alpha ((b_{s,t})) = (b_{{l_s}, {l_t}})^n_{s, t=1}
\]

\noindent where $\{l_1, \dots, l_n\}^2 \subseteq J$. If, additionally, $\{l_1,
\dots , l_n\}$ and $\{t_1, \dots , t_n\}$ are ordered subsets ($l_1 <
l_2$ etc), then we refer to $\gamma$ as an {\em ordered compression type
map}.

Let $A = A_1 \oplus \dots \oplus A_r$ where $A_1, \dots , A_r$ are the
matrix algebra summands of $A$ and let $\gamma : A \rightarrow B$ be a
map of compression type, as above, with $\gamma = \gamma_1 \oplus \dots
\oplus \gamma_p$, where $\gamma_1, \dots . \gamma_p$ are the elementary
summands of
$\gamma$. The map $\gamma$ has isometric restriction to $A_1$ if (and
only if) there is a summand $\gamma_i$ which is isometric on $A_1$. It
follows that $\gamma$ is isometric if (and only if) the summands can be
reordered and relabelled as $(\gamma_1 \oplus \dots \oplus \gamma_r)
\oplus (\gamma_{r+1} \oplus \dots \oplus \gamma_p)$ so that $\gamma^\pr =
\gamma_1 \oplus \dots \oplus \gamma_r$ is a (multiplicity one) algebra
injection of $A$ into $B$. With such a relabelling there is an
associated contractive left inverse map $\delta : B \rightarrow A$, of
compression type, satisfying $\delta \circ \gamma = id_A$. The map
$\delta$ is the composition of compression onto the range of $\gamma^\pr$,
followed by the inverse map of $\gamma^\pr$ restricted to its range.

Changing notation, let $A = \displaystyle{\lim_\to}\  (A_k, \phi_k),
$\ and $\  B = \displaystyle{\lim_\to}\  (B_k, \psi_k)$ be presentations of the
AF C*-algebras $A$ and $B$ where the maps $\phi_k$ and $\psi_k$ are
isometric C*-algebra injections. Assume that the matrix unit systems
$\{e^k_{ij} : (i,j) \in I_k\}$ and $\{f^k_{ij}: (i,j) \in J_k\}$ have
been chosen for $A_k$ and $B_k$ respectively, so that each map $\phi_k$
and $\psi_k$ maps matrix units to sums of matrix units, and, for
the moment, assume  that $A$ is unital and that the embeddings
 $\phi_k$ are unital.
Assume furthermore that the Bratteli diagram for the system $\{B_k, \psi_k\}$
has the Fermion property. By composing maps, forming a subsystem from
such compositions, and relabelling, we may assume that $B_k = M_{r_k}
\oplus B_k^\pr$ and that the partial embedding of $\psi_k$ from $M_{r_k}$
into $M_{r_{k+1}}$ has multiplicity at least two.
\vspace{.3in}

\newpage
\noindent {\bf Lemma 1.1.}\   \it With the assumptions above there is a
commuting diagram

\begin{center}
\setlength{\unitlength}{0.0125in}%
\begin{picture}(185,97)(120,685)
\thicklines
\put(240,760){\vector( 1, 0){ 40}}
\put(135,760){\makebox(0,0)[lb]{\raisebox{0pt}[0pt][0pt]{$A_1$}}}
\put(210,760){\makebox(0,0)[lb]{\raisebox{0pt}[0pt][0pt]{$A_2$}}}
\put(140,745){\vector( 0,-1){ 40}}
\put(220,745){\vector( 0,-1){ 40}}
\put(160,685){\vector( 1, 0){ 40}}
\put(240,685){\vector( 1, 0){ 40}}
\put(160,760){\vector( 1, 0){ 40}}
\put(305,720){\makebox(0,0)[lb]{\raisebox{0pt}[0pt][0pt]{\twlrm . . . }}}
\put(210,685){\makebox(0,0)[lb]{\raisebox{0pt}[0pt][0pt]{$B_{n_2}$}}}
\put(175,770){\makebox(0,0)[lb]{\raisebox{0pt}[0pt][0pt]{$\phi_1$}}}
\put(245,770){\makebox(0,0)[lb]{\raisebox{0pt}[0pt][0pt]{$\phi_2$}}}
\put(170,695){\makebox(0,0)[lb]{\raisebox{0pt}[0pt][0pt]{$\theta_1$}}}
\put(245,695){\makebox(0,0)[lb]{\raisebox{0pt}[0pt][0pt]{$\theta_2$}}}
\put(200,720){\makebox(0,0)[lb]{\raisebox{0pt}[0pt][0pt]{$\gamma_2$}}}
\put(120,720){\makebox(0,0)[lb]{\raisebox{0pt}[0pt][0pt]{$\gamma_1$}}}
\put(130,685){\makebox(0,0)[lb]{\raisebox{0pt}[0pt][0pt]{$B_{n_1}$}}}
\end{picture}
\end{center}

\noindent where each map $\gamma_k$  is an isometric linear map of compression
type relative to the given matrix unit system, and where $\theta_1,
\theta_2, \dots$ are compositions of the given embeddings $\psi_1,
\psi_2, \dots$. Furthermore there are linear contractions $\delta_k :
B_{n_k} \rightarrow A_k$, of compression type, satisfying $\delta_k
\circ \gamma_k = id$, such that the diagram

\begin{center}
\setlength{\unitlength}{0.0125in}%
\begin{picture}(185,97)(120,685)
\thicklines
\put(240,760){\vector( 1, 0){ 40}}
\put(135,760){\makebox(0,0)[lb]{\raisebox{0pt}[0pt][0pt]{$A_1$}}}
\put(210,760){\makebox(0,0)[lb]{\raisebox{0pt}[0pt][0pt]{$A_2$}}}
\put(160,685){\vector( 1, 0){ 40}}
\put(240,685){\vector( 1, 0){ 40}}
\put(140,705){\vector( 0, 1){ 40}}
\put(220,705){\vector( 0, 1){ 40}}
\put(160,760){\vector( 1, 0){ 40}}
\put(305,720){\makebox(0,0)[lb]{\raisebox{0pt}[0pt][0pt]{\twlrm . . . }}}
\put(210,685){\makebox(0,0)[lb]{\raisebox{0pt}[0pt][0pt]{$B_{n_2}$}}}
\put(175,770){\makebox(0,0)[lb]{\raisebox{0pt}[0pt][0pt]{$\phi_1$}}}
\put(245,770){\makebox(0,0)[lb]{\raisebox{0pt}[0pt][0pt]{$\phi_2$}}}
\put(170,695){\makebox(0,0)[lb]{\raisebox{0pt}[0pt][0pt]{$\theta_1$}}}
\put(245,695){\makebox(0,0)[lb]{\raisebox{0pt}[0pt][0pt]{$\theta_2$}}}
\put(200,720){\makebox(0,0)[lb]{\raisebox{0pt}[0pt][0pt]{$\gd_2$}}}
\put(120,720){\makebox(0,0)[lb]{\raisebox{0pt}[0pt][0pt]{$\gd_1$}}}
\put(130,685){\makebox(0,0)[lb]{\raisebox{0pt}[0pt][0pt]{$B_{n_1}$}}}
\end{picture}
\end{center}

\noindent commutes. In particular there exists an isometric injection $\gamma:
A
\rightarrow B$ and a contractive map $\delta: B \rightarrow A$ such that
$\gamma \circ \delta$ is a contractive projection onto the range of
$\gamma$.
\rm

\begin{mproof}   Using the Fermion property  for $B$ choose $n_1$
large enough so that there exists a multiplicity one linear isometry
$\gamma_1 : A_1 \rightarrow B_{n_1}$, of compression type, with range in
the summand $M_{r_1}$ of the decomposition $B_{n_1} = M_{r_1} \oplus
B^\prime_{n_1}$. We may assume that $\gamma_1$ has the form $\gamma_1 (a) =
[a \oplus 0_*] \oplus \{0\}$, that is that the partial embedding of $\gamma_1$
from $A_1$ into $M_{r_1}$ has a proper zero summand $0_*$.
(We indicate the distinguished first summand of $B_{n_2}$ with square brackets
and the remaining summands are grouped in braces.)
For $n_2 > n_1$, to
be chosen, the composed map $\theta_1: B_{n_1} \rightarrow B_{n_2}$ has
partial embeddings $\sigma_1: M_{r_1} \rightarrow M_{r_2}, \ \tau_1 :
B^\prime_{n_1} \rightarrow M_{r_2}, \ \tau_2 : B_{n_1} \rightarrow
B_{n_2}^\pr$, and by relabelling the matrix units of $M_{r_2}$ we may assume
that $\sigma_1$ is of standard type, that is,

\[
\sigma_1 (b) = b \oplus \dots \oplus b \op 0 \quad \quad
(\hbox{b appearing $t$ times}).
\]

\noindent (The zero summand may be absent.)
Using the Fermion property hypothesis we may choose $n_2$ so that $t$ is
arbitrarily large. Note that the composition $\theta_1 \circ \gamma_1$
has the form

\[
a \rightarrow [\sigma_1 (a \oplus 0_*) \oplus \tau_1 (0)] \oplus \{\tau_2
(\gamma_1 (a))\}.
\]

Consider now the given map $\phi_1 : A_1 \rightarrow A_2$. Let

\[
A_1 = A_{1,1} \oplus \dots \oplus A_{1,p}, \quad A_2 = A_{2,1} \oplus \dots
\oplus
A_{2,q}
\]

\noindent be the matrix algebra decompositions. Relabelling matrix units
of $A_2$ we may assume that $\phi_1$ is given in a standard form with
respect to the matrix unit systems, that is,

\[
\phi_1 : A_{1,1} \oplus \dots \oplus A_{1,p} \rightarrow A_{2,1}
\oplus \dots \oplus A_{2,q}
\]

\noindent where the summand of $\phi_1 (a_1 \oplus \dots \oplus a_p)$
in the matrix
summand $A_{2,t}$ is

\[
(\sum^{k_{1,t}}_1 \oplus a_1) \oplus \dots \oplus (\sum_1^{k_{p,t}}
\oplus a_p),
\]

\noindent with the understanding that some
of these summands may be absent. The
integer $k_{s,t}, $\ for $\ 1 \leq s \leq p, 1 \leq t \leq q,$ is the
multiplicity of the partial embedding for $\phi_1$ from $A_{1, s}$ to
$A_{2,t}$.

We now construct $\gamma_2 : A_2 \to B_{n_2}$ as an isometric
linear multiplicity one injection  of compression type, as
suggested by the following diagram:
\vspace{.3in}

\setlength{\unitlength}{0.0125in}%
\begin{picture}(150,107)(15,715)
\thicklines
\put( 85,720){\vector( 1, 0){ 35}}
\put( 50,785){\vector( 0,-1){ 45}}
\put(165,775){\vector( 0,-1){ 30}}
\put( 45,800){\makebox(0,0)[lb]{\raisebox{0pt}[0pt][0pt]{$a$}}}
\put( 15,715){\makebox(0,0)[lb]{\raisebox{0pt}[0pt][0pt]{$[a\op 0_*]\op 0$}}}
\put( 80,800){\vector( 1, 0){ 40}}
\put( 95,730){\makebox(0,0)[lb]{\raisebox{0pt}[0pt][0pt]{$\theta_1$}}}
\put(140,715){\makebox(0,0)[lb]{\raisebox{0pt}[0pt][0pt]
{$[\gs_1(a\op 0_*) \op \tau_1(0)] \op {\tau_2(\gamma_1(a))}$}}}
\put( 95,810){\makebox(0,0)[lb]{\raisebox{0pt}[0pt][0pt]{$\phi_1$}}}
\put( 25,760){\makebox(0,0)[lb]{\raisebox{0pt}[0pt][0pt]{$\gamma_1$}}}
\put(145,800){\makebox(0,0)[lb]{\raisebox{0pt}[0pt][0pt]
{$[\gS_{1,1}\op\dots\op\gS_{p,1}]\op\dots
\op[\gS_{1,q}\op\dots\op\gS_{p,q}]$}}}
\put(135,760){\makebox(0,0)[lb]{\raisebox{0pt}[0pt][0pt]{$\gamma_2$}}}
\end{picture}
\vspace{.3in}

\noindent Using the Fermion property choose $n_2$ large enough so that

\[
t \geq k_{i,1} + \dots + k_{i,q} \quad, 1 \leq i \leq p.
\]

\noindent  These
inequalities guarantee that there exists a one to one correspondence of
the summands of $\phi_1 (a)$ with {\em some} of the appropriate
nonzero summands of $\sigma_1 (a \oplus 0_*)$.
The point of this observation
is that there is a natural  multiplicity one algebra
injection $\gamma_2^\pr$ from $A_2$ to $M_{r_2}$ (of compression type)
which respects this correspondence.
Construct $\gamma_2$ by adding extra summands to $\gamma_2^\pr$  to obtain
a linear isometry, of compression type, satisfying
$\gamma_2(\phi_1(a)) = \theta_1(\gamma_1(a))$ for all $\ a$ in $\ A_1$.
Since $\phi_1$ is isometric this is possible.

Define $\delta_1 = \gamma_1^{-1} \circ \eta_1$ where $\eta_1$ is the
compression map onto the range of $\gamma_1$, noting that  $\gamma^{-1}_1$
is well-defined on this range. Similarly define $\delta_2$ in terms of
$\gamma_2$. Thus, $\delta_2 = (\gamma_2)^{-1} \circ \eta_2$ where $\eta_2$ is
the
compression onto the range of $\gamma_2$. To see the important equality
$\phi_1 \circ \delta_1 = \delta_2 \circ \theta_1$, let $b \in B_{n_1}$, and
let $\eta_1 (b) = [a \oplus 0_*] \oplus 0$. By the construction of
$\gamma_2$ and $\delta_2$ we have

\[
\delta_2 (\theta_1 (b)) = \delta_2 (\theta_1 (\eta_1 (b))).
\]

\noindent This is because the domain of $\delta_2$ is
{\em subordinate} to the nonzero summands
of $\sigma_1 (a \oplus 0_*)$, and $\theta_1(b -
\eta_1 (b))$ vanishes on these summands. Thus
\[
\delta_2 (\theta_1 (b))
= \delta_2(\theta_1(\eta_1(b))) = \delta_2 (\theta_1 ([a \oplus 0_*] \oplus 0))
= \phi_1 (a) = \phi_1
(\delta_1 (b)).
\]
Repeating the arguments above obtain  inductively  isometric maps\
$\gamma_3, \gamma_4, \dots$, with  distinguished contractive left inverses
$\delta_3, \delta_4, \dots$. Indeed, note that
after relabelling the matrix units of $B_{n_2}$
the map $\gamma_2 : A_2 \rightarrow
M_{r_2} \oplus B^\prime_{n_2}$ can be written in the form

\[
a \rightarrow [a \oplus 0 \oplus \gamma_{2,1} (a)] \oplus (\gamma_{2,2}
(a)),
\]

\noindent and so the construction of $\gamma_3^\pr, \delta_3, \gamma_3$\ is
obtained
in exactly the same way as $\gamma_2^\pr, \delta_2, \gamma_2$
\end{mproof}

The lemma shows that there is a complemented isometric linear injection
$A \rightarrow B$. If A is not unital then $A$ has a complemented
isometric linear injection into its unitisation $A^\sim$, and so the
nonunital case follows on consideration of the composition $A
\rightarrow A^\sim \rightarrow B$.

The next lemma also appears in Arazy's paper. For completeness we give
a proof. Write $ X \approx Y$\ if $X$ and $Y$ are linearly homeomorphic
Banach spaces.
\vspace{.3in}

\noindent {\bf Lemma 1.2.} \ \ $F \approx c_0(F)$.

\begin{mproof}  Realise $F$ as the direct limit
$\displaystyle{\lim_\to}\  (M_{2^k}, \rho_k)$ where $\rho_k : a \rightarrow
a \ot I_2$. Define a natural injection $\beta : c_0 (F)
\rightarrow F$ which is suggested by the following inclusion diagram.

\begin{center}
\setlength{\unitlength}{0.0125in}%
\begin{picture}(200,120)(80,680)
\thicklines
\put(160,740){\line( 1, 0){ 60}}
\put(220,740){\line( 0, 1){ 60}}
\put(220,740){\line( 0,-1){ 30}}
\put(220,710){\line( 1, 0){ 30}}
\put(250,710){\line( 0, 1){ 30}}
\put(250,740){\line(-1, 0){ 30}}
\put(250,710){\line( 1, 0){ 15}}
\put(265,710){\line( 0,-1){ 15}}
\put(265,695){\line(-1, 0){ 15}}
\put(250,695){\line( 0, 1){ 15}}
\multiput(265,695)(5.00000,-5.00000){2}{\makebox(0.4444,0.6667){\tenrm .}}
\put(160,800){\line( 1, 0){120}}
\put(280,800){\line( 0,-1){120}}
\put(280,680){\line(-1, 0){120}}
\put(160,680){\line( 0, 1){120}}
\put( 80,760){\makebox(0,0)[lb]{\raisebox{0pt}[0pt][0pt]{$F \ \ \ \supseteq $
}}}
\put(255,700){\makebox(0,0)[lb]{\raisebox{0pt}[0pt][0pt]{$.$}}}
\put(185,765){\makebox(0,0)[lb]{\raisebox{0pt}[0pt][0pt]{$F$}}}
\put(230,720){\makebox(0,0)[lb]{\raisebox{0pt}[0pt][0pt]{$F$}}}
\put(250,775){\makebox(0,0)[lb]{\raisebox{0pt}[0pt][0pt]{$0$}}}
\put(175,700){\makebox(0,0)[lb]{\raisebox{0pt}[0pt][0pt]{$0$}}}
\end{picture}
\end{center}

\noindent More precisely let $F_0$ be the subspace $\displaystyle{\lim_\to}\
(M_{2^k}^0
\rho_k)$, of codimension one given by the subsystem determined by the
subspaces
\[
M_{2^k}^0 = \{(a_{ij}) \in M_{2^k} :a_{2^k, 2^k} =0\}.
\]
Define
$c_0 (F) \rightarrow F_0$ as follows. Identify the first copy
of $F$ in $c_0 (F) = F \oplus F \oplus \dots$ with $p_1 F_0 p_1$
where $p_1 = e_{1,1}$ in $M_2$. (Identify $e_{1,1}$ with its image
in the limit.) Identify the second copy with $p_2 F_0
p_2$, where $p_2 = e_{33}$ in $M_{2^2}$, and so on. The resulting
inclusion $c_0(F) \rightarrow F$ has range which is the range of
the projection $E: F_0 \rightarrow F_0$ given by $E(a) = \lim_k
(\sum^k_{j=1} p_jap_j)$. Thus $c_0 (F)$ is complemented in $F_0$,
and hence in $F$. Thus $F \approx c_0 (F) \oplus X \approx c_0
(F) \oplus c_0 (F) \oplus X \approx c_0 (F) \oplus F \approx
c_0 (F)$.
\end{mproof}

The proof of the next theorem now reduces to a routine application of the
Pelczynski decomposition method.
\vspace{.3in}

\noindent {\bf Theorem 1.3.} \ \it Let $A$ and $B$ be AF C*-algebras given by
direct
systems whose Bratteli diagrams have the Fermion property. Then $A$ and
$B$ are isomorphic as topological vector spaces.
\rm

\begin{mproof} We may assume that $B = F$. By Lemma 1.1, and the remarks
concerning the unital case, there
exist contractive injective complemented maps $A \rightarrow F$ and $F
\rightarrow A$. Thus, by Lemma 1.2,

\[
c_0 (A) \rightarrow c_0 (F)
\approx F \rightarrow A.
\]

\noindent Hence, just as with $F$, we have $A \approx
c_0 (A) \oplus Y \approx c_0 (A) \oplus c_0 (A) \oplus Y
\approx c_0 (A) \oplus A \approx c_0 (A)$.

Consider now the fact that $A \approx F \oplus Z$ for some closed
subspace $Z$ of $A$, and obtain $A \approx F \oplus Z \approx c_0
(F) \oplus Z \approx c_0 (F) \oplus F \oplus Z \approx F \oplus A$.
Similarly, $F \approx F \oplus A$, and so $F \approx A$.
\end{mproof}

Let $A = \displaystyle{\lim_\to}\  A_k$ be an (infinite-dimensional)
AF C*-algebra, with a corresponding Bratteli diagram, which is simple.
This means that for each vertex $v$ of the diagram, at level $k$, there
is a lower level $m$ such that there exist downward paths from $v$
to {\em all} the vertices at level $m$. (See Bratteli \cite{bra}.)
In particular
any two vertices at a given level have downward paths that meet in a
common vertex. This weaker property is precisely the Bratteli
diagram criterion for the triviality of the centre of $A$. Suppose
additionally, that the Bratteli diagram fails to have the Fermion
property. Then there must exist a vertex with a {\em unique} downward
path. For otherwise there is repeated branching and convergence
characteristic of the Fermion property.

The unique downward path determines a subsystem of $A$ which defines a
subalgebra $J$ which is isomorphic to $\ck$ or $M_n$ for some $n$. Since
$J$ is in fact an ideal, and $A$ is simple, it follows  that $A= \ck$.
Thus we have obtained
\vspace{.3in}

\noindent {\bf Corollary 1.4.}\   \it Let $A$ be a simple
(infinite-dimensional)
approximately finite C*-algebra. Then $A \approx \ck$ or $A  \approx F$.
\rm
\vspace{.3in}

A C*-algebra is said to be of Type I if its star representations
generate Type I von Neumann algebras. Also it is known that this is
equivalent to the apparently weaker assertion that factorial star
representations are Type I.  (See, for example, \cite{fel-dor}.)
Using this we can
strengthen the last corollary.
\vspace{.3in}

\noindent {\bf Corrolary 1.5.}\   \it Let $A$ be an approximately finite
C*-alegbra which
is not  Type I. Then $A$ is isomorphic to $F$ as a linear topological
vector space.
\rm
\vspace{.3in}

\begin{mproof}  Let $A = \displaystyle{\lim_\to}\  A_k$, with Bratteli
diagram without the Fermion property. We show that the factorial
representations of $A$ are Type I. Note first that if $\pi : A
\rightarrow L (H)$ is a factorial representation, then $\ker \pi$ is an
ideal, and $A/\ker \pi$ is an AF C*-algebra with Bratteli diagram
obtained as a subdiagram of the diagram for $A$. (See Bratteli \cite{bra}.)
Since this
subdiagram also fails to have the Fermion property we may as well assume
that $\ker \pi = \{0\}$. With this assumption it follows that the centre
of $A$ must be trivial. By the argument preceding Corollary 4 it follows
that $A$ possesses an ideal which is isomorphic to $\ck$ and which is
associated with a vertex of the diagram that has a unique descending
path. Let $p$ be a minimal projection in the matrix summand
corresponding to this vertex. Then $p$ is minimal in $A$, and so $\pi
(p)$ is minimal in $\pi (A)^{\prime \prime}$. Thus the factor $\pi
(A)^{\prime \prime}$ is Type I. \hfill
\end{mproof}
\vspace{.3in}

\newpage
\noindent \bf Type I AF C*-algebras
\rm

\noindent
We comment on the isomorphism types of
the (infinite dimensional) Type I AF C*-algebras.

In the separable dual case the following three algebras present themselves:

(i) $\ c_0$, the space of diagonal compact operators,

(ii)$\ {\cal R} = (\sum_{k=1}^{\infty}\ \op M_k)_{c_0}$,

(iii) $\ck $, the compact operators.

To see that these are not isomorphic we can distinguish (i) from (ii) and (iii)
by noting that all bounded maps from $\ c_0$\ to $\ \ell^2$\ are
2-summing. A simple proof can be found in Pisier's notes \cite{pis-notes}.
On the other hand matrix realisations provide
 "top row" maps   $\ck \to \ell^2$,\ \  ${\cal R} \to \ell^2$
which are not 2-summing. Alternatively,
it can be seen in Hamana \cite{ham} and Chu and Iochum \cite{chu-ioc}
that (i) has the Dunford Pettis
property whereas (ii) and (iii) do not.
To distinguish ${\cal R}$ and $\ck$ one can note that the first space
has a dual space with the Schur property, that weakly convergent
sequences are norm convergent
whereas the trace class operators do not. (See \cite{ham} and \cite{chu-ioc}.)

If an AF C*-algebra has a separable dual space then it
is easy to see that it  does not have the
Fermion property and this limits the possibilities for the type of
Bratteli
diagram.
The Bratteli diagrams in this case fall naturally into three types,
namely, Type (i), in which there is a uniform bound on the
sizes of the matrix summands, Type (ii), in which for every path in the
diagram there is a uniform bound on the sizes of the associated summands,
and Type (iii), being the rest.
In the particular context of
a diagram with finite width
the associated algebra has a finite composition series and from this it
can be shown  to be isomorphic to $\ \ck$ as a Banach space.

Bessaga and Pelczynski \cite{bes-pel} have classified the spaces
$C(S)$ with $S$ countable. For each isomorphism type there is a countable
ordinal
$\alpha$ with $C(S) \approx C(\alpha)$, and $C(\beta) \approx C(\alpha)$
if and only if $\alpha \le \beta < \alpha^\omega$ or $\beta \le \alpha
< \beta^\omega$. In particular there is a continuum of isomorphism types
of abelian AF C*-algebras with diagrams of type (i).
One would expect there to be similar continua for the algebras whose diagrams
are of type (ii) and type (iii).

Turning to the algebras whose dual space is not separable
there are, in the first instance, six natural C*-algebras to consider.

(iv) $C(K), $ with $K$  a Cantor space,

(v) $C(K) \op {\cal R},$

(vi) $C(K) \ot {\cal R},$

(vii) $C(K) \op \ck$,

(viii) $(C(K) \ot  {\cal R}) \op  \ck$,

(ix) $ C(K) \ot \ck$.

These algebras are associated with six different
types of Bratteli diagram.
Thus, a Type (iv) diagram has
uncountably many paths and a uniform bound on the sizes of matrix summands.
A Type (v) diagram has uncountably many paths, is not of Type (iv),
has bounded matrix sizes on each path, and for each integer $n$
has only countably many paths on
which the matrix sizes exceed any $n$.
Similar descriptions hold for the remaining types, culminating in Type (ix)
for which there are uncountably many paths each with unbounded matrix sizes.
Moreover each Bratteli diagram without the Fermion property, which is not of
Types (i), (ii) or (iii) is one of these six types.
It seems plausible that isomorphic AF C*-algebras have Bratteli diagrams
of the same diagram type, and within some of these diagram types one would
expect there to be added ordinal type complexities as in the
separable dual case. (For example one might expect there to be a
continua
of Banach space types of the form $C(K) \op (C(\alpha) \ot \ck)$).

With regard to the duals of AF C*-algebras, Wojtaszczyk has shown that
there are just three separable duals, namely the duals of the algebras
(i), (ii), and (ii). It seems reasonably to conjecture that there
are precisely nine dual spaces of Type I AF C*-algebras.

\section{\bf Triangular subalgebras}

\noindent There are three well-known families of triangular subalgebras of UHF
C*-algebras, namely the refinement algebras
$ \displaystyle{\lim_\to}\  (T_{n_k},\rho_k)$, the standard
algebras $ \displaystyle{\lim_\to}\  (T_{n_k},\sigma_k)$
and the
alternation algebras,
$ \displaystyle{\lim_\to}\  (T_{n_k},\alpha_k)$.  The
(unital) embeddings determining these limits have the form

\[
\rho_k ((a_{ij})) = (a_{ij} I_{t_k}),\ \
\sigma_k (a) = I_{t_k} \otimes a,
\]

\noindent where $t_k$ is the multiplicity of the embedding, and in the
alternation
case $\alpha_k$ alternates between these two types. We shall prove the
following theorem.
\vspace{.3in}

\noindent {\bf Theorem 2.1.}\    \it (i)The standard limit
algebras are isomorphic as Banach spaces.
(ii) The refinement limit algebras are isomorphic as Banach
spaces.
(iii) The alternation limit algebras are isomorphic as Banach
spaces.
\rm
\vspace{.3in}

Another well-known class consists of the various ``refinement with
twist'' limits $\displaystyle{\lim_\to} (T_{n_k}, \tau_k)$,
where $\tau_k$ agrees
with $\rho_k$ on all the standard matrix units $e_{ij}$ of $T_{n_k}$,
with the exception of those superdiagonal matrix units in the last
column. For these

\[
\tau_k (e_{i, n_k}) = e_{i,n_k} \otimes u_k
\]

\noindent where $u_k$ is a permutation unitary in $M_{t_k}$.
It was shown in Hopenwasser and Power \cite{hop-pow} that
these algebras provide uncountably
many algebra isomorphism classes, distinct from the
refinement limits. On the other hand we have
\vspace{.3in}

\noindent {\bf Theorem 2.2.}\  \it If $A$ is a refinement
with twist algebra, as above,
then as a Banach space, $A$ is isomorphic to the model refinement algebra
$\ct_{2^\infty}$.
\rm
\vspace{.3in}

The algebras above are examples of (canonical regular) triangular
subalgebras of UHF C*-algebras. In a different direction
one can generalize the
standard embedding limit algebras by considering ordered Bratteli
diagrams (\cite{scp-k0}, \cite{poo-wag}.) A typical such embedding has the form

\[
\beta : T_{q_1} \oplus \dots \oplus T_{q_s} \rightarrow T_{p_1} \oplus
\dots \oplus T_{p_r}
\]
where
\[
\beta : a_1 \oplus \dots \oplus a_s \rightarrow (\sum^{t_1}_{j=1} \oplus
b_{1,j} ) \oplus \dots \oplus (\sum^{t_r}_{j=1} \op b_{r,j})
\]

\noindent and where each $b_{k,l}$ is one of the summands $a_1, \dots, a_s$.
The
summations here mean block diagonal direct sums. In the nonunital case
one also allows the $b_{k,l}$ to be zero summands. For a simple example,
consider the embedding $\beta_1$ from $T_2 \oplus T_3 \oplus T_4$ to
$T_7 \oplus T_6 \oplus T_5$ given by
\[
\beta_1 : a \oplus b \oplus c \rightarrow (b \oplus c) \oplus (a \oplus
c) \oplus ( a \oplus b).
\]

\noindent This embedding is not inner conjugate to
\[
\beta_2 : a \oplus b \oplus c \rightarrow (b \oplus c) \oplus (c \oplus
a) \oplus (a \oplus b)
\]
and the difference can be indicated by ordered Bratteli diagrams. A consequence
of this diversity is the existence of uncountably many nonisomorphic
limit algebras with the same generated C*-algebra. Once again, however,
they correspond to a unique Banach space type.

\vspace{.3in}

\noindent {\bf Theorem 2.3.}\  \it Let $A$ be a triangular limit algebra
determined by an
ordered Bratteli diagram which has the Fermion property. Then, as a
Banach space, $A$ is isomorphic to the model algebra\ \
$\cs_{2^\infty}$.
\rm
\vspace{.3in}

Recall the definition of a linear map $\gamma : M_n \rightarrow M_m$
which is of {\it ordered compression type} and note that such a map has
a restriction $\gamma : T_n \rightarrow T_m$. Using direct sums of such
maps define {\it ordered compression type maps} $T_q \rightarrow T_{p_1}
\oplus \dots \oplus T_{p_r}$ and
use these to define general {\it ordered compression type
maps} $\gamma = T_{q_1} \oplus \dots \oplus T_{q_s} \rightarrow T_{p_1}
\oplus \dots \oplus T_{p_r}$. If $\gamma (a_1 \oplus \dots \oplus a_s)$
has at least one complete summand $a_i$, for each $i$, then $\gamma$ is
isometric. In this case it follows, as in section 2, that $\gamma$ has an
associated left inverse.
\vspace{.3in}

\noindent {\it The Proof of Theorem 2.3.} \ \ Note first that the conclusions
of Lemma 1.1 hold in
the triangular context wherein we make the following new assumptions :
$A = \displaystyle{\lim_\to}(A_k, \varphi_k)$ and $B =
\displaystyle{\lim_\to}(B_k,
\psi_k)$ are limit algebras determined by ordered Bratteli diagrams and
the diagram for $B$ has the Fermion property. To adapt the proof
we use ordered compression maps in place of compression maps.
The argument is virtually the same but notationally awkward
since we cannot make simplifying reorderings of summands by relabelling
matrix units in the codomain. The first occasion for this is the
expression for $\theta_1 \circ \gamma_1 (a)$. (The  map $\gamma_1$ can
be chosen to be of ordered compression type.)
However $\theta_1 \circ \gamma_1 (a)$ is an ordered direct sum such
that, in the first summand $T_{r_2}$ of $B_{n_2}$, for suitably large
$n_2$, there appear many summands which are copies of the summands of
$a$. Thus, as before, there is an association of the summands of
$\phi_1 (a)$ with some of the summands of $\theta_1 \circ \gamma_1
(a)$, if $n_1$ is large enough. However we assume additionally that this
association respects the order in which summands of $\phi_1 (a)$ appear
in each summand $A_{2,i}$. As before define a multiplicity one injection
$\gamma_2^\pr$\ of ordered compression type, from $A_2$ to $B_{n_2}$, which
respects
this correspondence. The extension $\gamma_2$ of $\gamma_2^\pr$
and the left inverse $\delta_2$ of $\gamma_2$ (and $\gamma_2^\pr$) are
defined as before, and once again the desired commuting diagrams follow
upon iterating this procedure. In the case that $A$ is nonunital obtain
a complemented unital isometric injection $A \rightarrow B$ by
considering the unitisation of $A$.

The argument of Lemma 1.2 also serves to give a natural complemented
injection $c_0 (\cs_{2^\infty}) \rightarrow \cs_{2^\infty}$.
(The diagram of Lemma 1.2,
however, is not appropriate for standard embeddings.) The remainder of
the proof follows as before. \hfill $\Box$
\vspace{.3in}

\noindent {\it The proof of Theorem 2.1 (ii)}.
\ \ This follows a similar scheme. Let $A = \displaystyle{\lim_\to}(A_k,
\phi_k)$
be a refinement limit algebra and
let $\ct_{2^\infty} = \displaystyle{\lim_\to}(B_k,
\psi_k)$ be the $2^\infty$ refinement limit algebra. To obtain the
appropriate version of Lemma 1.1 first choose $n_1$ large enough so that
there is a multiplicity one ordered compression type map $\gamma_1 : A_1
\rightarrow B_{n_1}$ given by $\gamma_1 (a) = a \oplus 0$ (block
diagonal direct sum). The matrix $\phi_1 (a)$ has the form $(a_{ij}
I_t)$ where $t$ is the multiplicity of $\phi_1$. Choose $\theta_1 :
B_{n_1} \rightarrow B_{n_2}$, a composition of consecutive maps $\psi_1,
\psi_2, \dots$ so that the multiplicity of $\theta_1$ exceeds that of
$\phi_1$. Thus

\[
\phi_1 ((a_{ij})) = (a_{ij} I_t)
\]
and
\[
\theta_1 (\gamma_1 ((a_{ij}))) = (a_{ij} I_s) \oplus 0
\]

\noindent where $s > t$. There is now natural multiplicity one isometric
(algebra)
injection $\gamma^\prime_2 : A_2 \rightarrow B_{n_2}$ with the property
that

\[
\gamma^\prime_2 (a_{ij} I_t) = (a_{ij} (I_t \oplus 0_{s-t})) \oplus 0.
\]

\noindent This (as before) does not yet give a commuting square. Nevertheless
we
can add multiplicity one summands of ordered compression type to create
$\gamma_2$, an extension of $\gamma_2^\prime$, so that $\gamma (a_{ij}
I_t) = (a_{ij} I_s) \oplus 0$, and thus we obtain the first commuting
square of Lemma 1.1 in this context. Furthermore, $\gamma_1$ has a natural
left inverse $\delta_1$, and, as before, $\gamma_2^\prime$ can be used
in the definition of a left inverse $\delta_2$ for $\gamma_2$, which
extends $\delta_1$ in the obvious way. The construction of the desired
maps $\gamma_3, \delta_3, \gamma_4, \delta_4, \dots$ is obtained
similarly.

We have $\ct_{2^\infty} \approx c_0 (\ct_{2^\infty})$, by the
argument of Lemma 1.2 (the
diagram is appropriate this time), and the proof is completed as
before. \hfill $\Box$
\vspace{.3in}

\noindent {\it The proof of Theorem 2.1 (iii)}
Let $S$ and $r$ be the generalised integers associated
with the triangular limit algebras $\cs_s$ and $\ct_r$.
\ \ Consider the subalgebra $\cs_s \star \ct_r$ of $B_s \ot B_r$
given by

\[
\cs_s \star \ct_r = (\cs_s \cap \cs_s^*) \ot \ct_r \ + \  (\cs_s^0 \ot
B_r)
\]

\noindent where $B_s = C^*(\cs_s)$ and $B_r = C^*(\ct_r)$
are the generated C*-algebras, and where $\cs_s^0$ is the strictly upper
triangular subalgebra of $\cs_s$. In the terminology of
\cite{scp-book} and \cite{scp-lexico} this algebra is the lexicographic
product of the ordered pair $\cs_s, \ct_r$.
If $s$ and $r$ are not finite then this product coincides
with the proper alternation algebra
for the pair $r, s$.
This formula forms the basis of the proof.

Note first that $\cs_s \star \ct_r$ is bicontinuously isomorphic to the
$\ell^\infty$ direct sum of the component spaces above. Also $\cs_s \cap
\cs_s^*$ is
C*-algebraically isomorphic to $C(K)$, where $K$ is a Cantor space.
{}From the linear isomorphism
$\ct_r \approx \ct_{2^\infty}$ we obtain

\[
(\cs_s \cap \cs_s^*) \ot \ct_r \approx
C(K) \ot \ct_r \approx C(K) \ot \ct_{2^\infty}
\approx (\cs_s \cap \cs_s^*) \ot \ct_{2^\infty}
\]

\noindent It remains then to show that $\cs_r^0 \ot B_r$
and $\cs_{2^\infty} \ot B_{2^\infty}$
are bicontinuously isomorphic. However the maps of Lemma 1.1 and its
non-self-adjoint variants respect tensor products.
For example, the isometric map $\gamma : A \to B$ of
Lemma 1.1 also provide complemented isometric injections
$\gamma \ot id : A \ot D \to B \ot D$, where
$D$ is a closed subspace of an AF C*-algebra (for example)
and the tensor product is the injective, or spatial, tensor product.
And so we obtain the needed equivalences

\[
\cs_r^0 \ot B_r \approx \cs_{2^\infty}^0 \ot B_r \approx \cs_{2^\infty}^0
\ot B_{2^\infty}
\]

\noindent and the proof follows.
\hfill $\Box$.

\vspace{.3in}

\noindent \it The proof of Theorem 2.2. \rm
\ \ Let
$ A = \displaystyle{\lim_\to}(T_{n_k},\tau_k)$ be a refinement with twist
limit algebra.
Let $A_0 = \displaystyle{\lim_\to}(T_{n_k}^-,\tau_k)
= \displaystyle{\lim_\to}(T_{n_k}^-,\rho_k)$
where $ T_{n_k}^- = \mbox{span}\{e_{i,j}: j < N_k\}$.
Also, let $\ct = \displaystyle{\lim_\to}(T_{n_k},\rho_k)$ be the associated
refinement algebra.
Then $T_{n_k} \approx T_{n_k}^- \op \ell_{n_k}^{2}$.
That is, the map

\[
a \to a(1 - e_{n_k,n_k}) \op ae_{n_k,n_k}
\]

\noindent is a bicontinuous linear isomorphism to the $\ell^{\infty}$ direct
sum.
{}From the commuting diagram

\begin{center}
\setlength{\unitlength}{0.0125in}%
\begin{picture}(185,97)(120,685)
\thicklines
\put(240,760){\vector( 1, 0){ 40}}
\put(135,760){\makebox(0,0)[lb]{\raisebox{0pt}[0pt][0pt]{$T_{n_1}$}}}
\put(210,760){\makebox(0,0)[lb]{\raisebox{0pt}[0pt][0pt]{$T_{n_2}$}}}
\put(140,745){\vector( 0,-1){ 40}}
\put(220,745){\vector( 0,-1){ 40}}
\put(165,685){\vector( 1, 0){ 35}}
\put(245,685){\vector( 1, 0){ 35}}
\put(160,760){\vector( 1, 0){ 40}}
\put(305,720){\makebox(0,0)[lb]{\raisebox{0pt}[0pt][0pt]{\twlrm . . . }}}
\put(205,685){\makebox(0,0)[lb]
{\raisebox{0pt}[0pt][0pt]{$T_{n_2}^- \op {\ell_{n_2}^{2}}$}}}
\put(175,770){\makebox(0,0)[lb]
{\raisebox{0pt}[0pt][0pt]{$\tau_1$}}}
\put(245,770){\makebox(0,0)[lb]{\raisebox{0pt}[0pt][0pt]{$\tau_2$}}}
\put(125,685){\makebox(0,0)[lb]
{\raisebox{0pt}[0pt][0pt]{$T_{n_1}^- \op {\ell_{n_1}^{2}}$}}}
\end{picture}
\end{center}

\noindent we obtain a bicontinuous isomorphism $A \approx A_0 \op \ell^2$.
Similarly $\ct \approx A_0 \op \ell^2$
and so the theorem now follows from Theorem 2.1.
\hfill $\Box$.
\vspace{.3in}

\noindent {\bf Remark 2.4.} \  With respect to injections which map
matrix units to sums of matrix units the Banach space $\cs_{2^\infty}$
is not as injective as $\ct_{2^\infty}$. More specifically if
$\phi :M_2 \to T_{m_1}$ and
$\phi^\pr :M_n \to T_{m_2}$ are isometric linear
maps which map matrix units to sums
of matrix units, and if $i : M_2 \to M_n$ is a unital C*-algebra injection,
and if $\phi^\pr \circ i = \sigma \circ \phi$, then
it can be shown that $2n \le m_1$. This suggests that there are no complemented
injections of the Fermion algebra in $\cs_{2^\infty}$,
whereas it is clear that there
are such injections for $\ct_{2^\infty}$.

\end{document}